\documentstyle[12pt]{report}

\renewcommand\a{\alpha}
\renewcommand\b{\beta}

\renewcommand\d{\dots}
\newcommand\dd{\diamondsuit}

\newcommand\e{\eta}

\newcommand\g{\gamma}
\newcommand\G{{\cal G}}

\renewcommand\H{{\cal H}}

\newcommand\J{{\cal J}}
\newcommand\k{\kappa}

\newcommand\la{\langle}
\renewcommand\ll{\lambda}
\newcommand\m{\mu}
\newcommand\N{{\cal N}}
\renewcommand\O{\Omega}
\newcommand\oo{\omega}

\newcommand\op[1]{\mathop{\rm #1}\nolimits}
\newcommand\ot{\otimes}
\newcommand\p{\partial}
\newcommand\po{$\!\!\!{\bf .}$ }

\newcommand\R{{\rm I\hspace{-2.5pt} R}}
\newcommand\ra{\rangle}
\renewcommand\t{\times}
\newcommand\te{\theta}

\newcommand\w{\wedge}
\newcommand\we{\wedge}
\newcommand\x{\xi}

\newcommand\Z{{\rm Z\mkern-5muZ}}

\def\Rom#1{\uppercase\expandafter{\romannumeral#1}}

\newcommand\1{{\bf 1}}

\newcommand\qed{\phantom{\underline{y}}\hfill\hfill$\Box$}

\newcommand\bib[1]{\bibitem[#1]{#1}}

\newcommand{\text}[1]{{\mbox{\rm #1}}}
\newcommand{\dfrac}[2]{\frac{\displaystyle #1}{\displaystyle #2}}

\newtheorem{th}{Theorem}
\newtheorem{prop}{Proposition}
\newtheorem{lem}{Lemma}

{\endtrivlist}
\newenvironment{cor}{\trivlist \item[\hskip \labelsep{\bf Corollary.}]\it}%
{\endtrivlist}
\newenvironment{dfn}[1]{\trivlist \item[\hskip \labelsep{\bf Definition #1.}]}%
{\endtrivlist}
{\endtrivlist}
\newenvironment{rk}{\trivlist \item[\hskip
\labelsep{{\it\underline{Remark}.\/}}]}%
{\endtrivlist}
\newenvironment{proof}{\trivlist \item[\hskip
\labelsep{{\it\underline{Proof}.\/}}]}%
{\endtrivlist}
{\endtrivlist}

\newcounter{a}
\newcounter{f}
\makeatletter
\newcommand{\@thefnmark}{$^\fnsymbol{f}$}
\renewcommand{\@makefnmark}{\hbox{\mathsurround=0pt
                           $^{\fnsymbol{f}}$}}
\renewcommand{\@makefntext}[1]{\parindent=1em\noindent
            \hbox to 1.8em{\hss$^{\fnsymbol{f}}$}#1}
\makeatother

\begin{document}

\title{Symplectic and contact Lie algebras with an application to
Monge-Amp\'ere equations}

\author{\bf Boris Kruglikov}

\date{August 10, 1997}

\maketitle

\begin{abstract}
In this paper we consider symplectic and contact Lie algebras.
We define contactization and symplectization procedures and
describe its main properties. We also give classification
of such algebras in dimensions 3 and 4. The classification in
dimension~4 is closely connected with normal forms of nondegenerate
elliptic equations of the second order on two-dimensional surfaces
with transitive symmetry group in first jets. We point out this
connection and discuss normal forms.
\end{abstract}

\tableofcontents

\clearpage

%%%%%%%%%%%%%%%%%%%%%%%%%%%%%%%%%%%%%%%%%%%%%%%%%%%%%%%%%%%%%%%%%%%%%%%%%%

\chapter*{Introduction}
\addcontentsline{toc}{chapter}{\bf\quad \  Introduction}

\hspace{13.5pt}
A Lie-Poisson group $G$ is a Lie group equipped with a Poisson structure
such that if we extend the Poisson structure naturally to the product $G\t
G$ then the Lie multiplication map $G\t G\to G$ is Poisson. Since the
structures are invariant we may equivalently talk of Lie-Poisson algebras.
When the Poisson structure is nondegenerate we call it symplectic. Thus we
obtain the problem of description of symplectic Lie algebras.

In this paper we define and study symplectic and contact Lie algebras. For
general manifolds in~\cite{A} the procedure of symplectization and
contactization was defined which is a functor between the categories of
exact symplectic and contact manifolds. We study the corresponding notions
for contact and exact symplectic Lie algebras. Note that the invariance
condition forces us to change the general constructions. We also give
a description of contact in dimension~3 and symplectic in dimension~4 Lie
algebras.

Recall (\cite{D}) that a Lie bialgebra is a Lie algebra $\G$ with a
Lie algebra structure on $\G^*$, these structures being compatible.
There is a bijective correspondence between Lie-Poisson algebras and Lie
bialgebras. When Lie-Poisson algebra is nondegenerate, i.e. symplectic, we
obtain that Lie algebra structures on $\G$ and $\G^*$ are equivalent or
that the map $\G\to\G\w\G\subset\G\ot\G$ dual to the Lie multiplication on
$\G^*$ is a 1-cocycle. Thus every symplectic Lie algebra gives rise to a
solution of the classical Yang-Baxter equation (\cite{D}).

Another application of symplectic Lie algebras occurs in dimension~4 where
they correspond to invariant under some Lie group transitive action
elliptic equations depending on two variables. To be more exact we talk
about Monge-Amp\'ere equations on the plane and about the generalization of
them given in~\cite{L}; we called the corresponding equations generalized
Monge-Amp\'ere. For such equations the equivalence problem was solved
in~\cite{K}. Namely for every nondegenerate elliptic equation there was
constructed some canonical $\{e\}$-structure and the defining Monge-Amp\'ere
equation structures proved to have a canonical expression by means of this
$\{e\}$-structure. Thus the recovering procedure for Monge-Amp\'ere equation
by its invariant is close to algorithmic.
Basing on the classification of 4-dimensional symplectic Lie
algebras we may give normal forms of nondegenerate elliptic generalized
Monge-Amp\'ere equations.

The author is grateful to professor V.\,V.~Lychagin for a warm attention to
the work and helpful discussions.%
\footnote{
The research described in this publication was made possible in part by
Grant RFFI $\N$ 96-010710 of Russian Fund of Fundumental Researches.}

%%%%%%%%%%%%%%%%%%%%%%%%%%%%%%%%%%%%%%%%%%%%%%%%%%%%%%%%%%%%%%%%%%%%%%%%%%
% 1 %
 \chapter{Symplectic and contact Lie algebras}

\hspace{13.5pt}
It is natural to consider left-invariant structures on Lie groups.
For in this case they may be treated as structures on Lie algebras. For
example the de Rham's complex for left-invariant differential forms on a
Lie group $G$ coincides with the complex for Lie algebra $\G$ cohomologies'
complex with coefficients in the trivial module $\R$. Hence everything may
be expressed in the language of algebras (and re-presented for groups).

%%%%%%%%%%%%%%%%%%%%%%%%%%%%%%%%%%%%%%%%%%%%%%%%%%%%%%%%%%%%%%%%%%%%%%%%%%
% 1.1 %
\section{The definitions and examples of
nonexistence}

\hspace{13.5pt}
Let us consider a Lie algebra $\G$ and the corresponding cohomology complex
with differentials $d: C^i\to C^{i+1}$ (see~\cite{J}, \cite{SL}),
$C^i=\w^i\G^*$.

 \begin{dfn}{1}
Lie algebra $\G$ of even dimension $2n$ is called symplectic when equipped
with a nondegenerate closed 2-form $\oo\in C^2$: $d\oo=0$, $\oo^n\ne0\in
C^{2n}$. Lie algebra $\G$ of odd dimension $2n+1$ is called contact when
equipped with a 1-form $\a\in C^1=\G^*$ such that
$\a\w(d\a)^n\ne0\in C^{2n+1}$.
 \end{dfn}

Let us also give an equivalent definition (the only difference for contact
structures is that $\Pi$ determines the projective class
$\{s\a\,\vert\,s\in\R^*\}$).

 \begin{dfn}{$1'$}
A closed 2-form $\oo$ is called symplectic if the mapping $\G\to\G^*$,
$v\mapsto i_v\oo$, is an isomorphism (of linear spaces). A Lie algebra is
called contact if it possesses a codimension 1 linear subspace $\Pi$ such
that for none vector $w\in\Pi\setminus\{0\}$ it holds $ad_w\Pi\subset\Pi$.
 \end{dfn}

Now the description of Lie groups with left-invariant symplectic (contact)
structure is reduced to the description of symplectic (contact) Lie
algebras. Not every even-dimensional Lie algebra is symplectic (the same
for contact of course, take the commutative one).
For example let $\G$ be a compact Lie
algebra (whence reductive, see~\cite{VO}). Then it cannot possess an exact
symplectic form, for this is obvious for (compact) Lie groups in view of
Stokes' formula. Thus for compact algebras with $H^2(\G)=0$ we derive they
are not symplectic. Thanks to the K\"unneth's formula this is the case with
compact semisimple or compact reductive with 1-dimensional center Lie
algebras. Here we used the Whitehead lemmas (\cite{J}): for semisimple Lie
algebra $\H$ we have $H^1(\H;\R)=H^2(\H;\R)=0$.

Now the compactness restriction above may be easily got over. Let us remind
that a Lie algebra $\G^n$ is called unitary (\cite{F}) if $H_n(\G;\R)\ne0$.
Equivalently for some (=every) basis $e_i$ of $\G$ and the corresponding
basis $e_i^*$ of $\G^*$ we have $\sum de_i^*(e_i,\cdot)=0$. As examples we
have reductive and nilpotent Lie algebras. For unitary algebras the
Poincar\'e duality holds: $H^k(\G)\simeq H_{n-k}(\G)$.

 \begin{th}\po
Unitary Lie algebra $\G$ possesses no exact symplectic form.
 \end{th}

 \begin{proof}
If the symplectic form $\oo$ on $\G$ with $\op{dim}\G=2n$ is exact,
$\oo=d\a$, then the differential of the form $\a\w\oo^{n-1}$ is nontrivial
and hence $H^{2n}(\G)=0$, which contradicts the Poincar\'e duality
since $H_0(\G)=\R$. \qed
 \end{proof}

 \begin{cor}
If for a unitary Lie algebra $\G$ the second cohomology group is trivial,
$H^2(\G)=0$, then it is not symplectic. These are the cases with semisimple
and reductive with 1-dimensional center Lie algebras. \qed
 \end{cor}

%%%%%%%%%%%%%%%%%%%%%%%%%%%%%%%%%%%%%%%%%%%%%%%%%%%%%%%%%%%%%%%%%%%%%%%%%%
% 1.2 %
\section{Unification of symplectization
and contactization}

\hspace{13.5pt}
The following construction coincides with symplectization for even
dimension and contactization for odd.

 \begin{dfn}{2}
The genre of a differential form $\a\in\O^1(M)$ is the highest degree of a
nonzero form in the sequence
 $$
\a, d\a, \a\w d\a, (d\a)^2, \a\w(d\a)^2, \d
 $$
Let us call a pair $(M^n,\a)$ nondegenerate if the genre of $\a$ equals
$n=\op{dim}M$.
 \end{dfn}

For even $n$ we get an exact symplectic manifold with the structure $d\a$
and for odd we get a contact one with the structure form $\a$
(to be exact a contact structure is a distribution $\op{Ker}\a$,
but up to a double covering any contact
structure may be obtained from a contact form and keeping in mind
Lie algebra constructions we see that these two ways are equivalent).

Now we define {\it the suspension\/}, i.e. a method to construct
an $(n+1)$-dimensional nondegenerate manifold $(M_+,\a_+)$
by a nondegenerate $(M,\a)$. It is called symplectization in
odd case and contactization in the even (see~\cite{A}). Let $M_+^{n+1}= M^n
\t\R^*$, where $\R^*=\R\setminus\{0\}=\R_-\cup\R_+$ (one may take instead
$\R_+$). Let $t$ be a parameter on $\R^*$ and define the form $\a_+$ to be
$t\a+dt$, where we consider $\a$ to be naturally lifted to $M_+$.

 \begin{prop}\po
The pair $(M_+,\a_+)$ is nondegenerate.
 \end{prop}

 \begin{proof}
$d\a_+=dt\w\a+td\a$. Hence for $n=2k$ we have $\a_+\w(d\a_+)^n= t^ndt\w
(d\a)^n\ne0$. And for $n=2k+1$ $(d\a_+)^{n+1}= t^ndt\w\a\w (d\a)^n\ne0$.
\qed
 \end{proof}

Another way to obtain this construction is to consider $M_+$ as
a 1-dimensio\-nal bundle over $M$ with the connection form $\a_+$ and the
curvature form $d\a$.

Now we may define the prolongations of the structure diffeomorphisms
(contactizations of symplectomorphisms and symplectizations of
contactomorphisms) and so on.

%%%%%%%%%%%%%%%%%%%%%%%%%%%%%%%%%%%%%%%%%%%%%%%%%%%%%%%%%%%%%%%%%%%%%%%%%%
% 1.3 %
\section{Symplectizations and contactizations of
Lie algebras}

\hspace{13.5pt}
The differential constructions considered above are not invariant and hence
must be a little changed for Lie algebras. It turns out that every time
there exists some contactization while symplectization sometimes cannot take
place. By the direct sum beyond we mean the sum as vector spaces not algebras.
Let us denote by $Z^i$ the subspace $\op{Ker}d$ of $C^i$.

Let us recall that Lie algebra structure can be fixed in two equivalent ways:
by the commuting relations $[e_i,e_j]=c_{ij}^ke_k$ in Lie algebra or
by the Maurer-Cartan equations $df^k=-\dfrac12c_{ij}^kf^i\w f^j$ in its
dual, where $f^i$ stands for the dual basis to $e_i$. In addition the
Jacobi identities are equivalent to the integrability conditions $d^2=0$.
\newline
%%%%%%%%%%%%%%%%%%%%%%%%%%%%%%%%%%%%%%%
\newline
{\large\bf 1.3.1.\  \ Suspensions}
\bigskip

Let $\G^{n+1}=\H^n\oplus\R^1$, $\H$ being an ideal (it is also possible to
consider the case of subalgebras). Let us suppose that $(\H,\a)$ is a
nondegenerate pair. Set $\a_+=\pi^*\a$ for the projection $\pi:\G\to\H$
induced from the direct decomposition.

 \begin{dfn}{3}
The pair $(\G,\a_+)$ is called {\it a suspension\/} over $(\H,\a)$ if it is
a nondegenerate pair. For the even-dimensional case it is also called
a contactization, for the odd one --- a symplectization.
 \end{dfn}

Since $\H$ is an ideal the Lie algebra structure on $\G$ is given by an
outer derivation $[A]\in H^1(\H,\H)$. Here $[A]$ is a class of the
derivation $A\in\op{Der}(\H)$ modulo the subalgebra of inner
derivations $\{ad_w\,\vert\,w\in\H\}$ and the derivation $A$
stands for $ad_v$ with $v$ a generator of $\R^1$.

Note that under the change of the transversal line $\R^1$, $v\mapsto v+w$,
$w\in\H$, the suspension form is changed by the formula
$\a_+\mapsto\a_+-\a(w)v^*$, where $v^*\in(\R^1)^*$ is the dual to $v$
covector, $v^*(\H)=0$, $v^*(v)=1$. Thus it seems that the set of parameters
for the suspension is the derivations $A\in\op{Der}(\H)$ themselves, not
their classes.

Let us take now a geometrical point of view.
\newline
{\it\underline{Contactization}.\/}
Let $(\H^{2n},d\a)$ be a symplectic Lie algebra, $\a\in\H^*$.
Let $\Pi^{2n-1}=\op{Ker}\a$, $u$ be some transversal vector to it in $\H$
and $w\in\op{Ker}(\left.d\a\right\vert_\Pi)\setminus\{0\}$. The condition
$\a\w(d\a)^n\ne0$ on $\G$ is equivalent to $d\a(v,w)\ne0$ or
$ad_vw\notin\Pi$. Note that the addition to $A=ad_v$ an inner
derivation $ad_x$, $x\in\Pi$, does not change the form $\a_+$.
Moreover the derivation $A=ad_u$ enjoys the property $Aw\notin\Pi$.
Thus we have proved

 \begin{prop}\po
The set of all contactizations is nonempty and is para\-met\-ri\-zed by the set
 $$
\,\quad\qquad\qquad\qquad\hfill\hfill
\{A\in\op{Der}(\H)\,\vert\,Aw\notin\Pi\} / \{ad_x\,\vert\,x\in\Pi\}.
\,\quad\qquad\qquad\qquad\hfill\hfill\Box
 $$
 \end{prop}

Note that the contactization may also be applied for any other potential
$\a+\b$, $d\b=0$, of the symplectic form $\oo=d\a$.
\newline
{\it\underline{Symplectization}.\/}
Let $(\H^{2n-1},\a)$ be a contact Lie algebra, $\a\in\H^*$.
Let $\Pi^{2n-2}=\op{Ker}\a$ and $w\in\op{Ker}d\a\setminus\{0\}$.
The condition $(d\a)^n\ne0$ on $\G$ is equivalent to $d\a(v,w)\ne0$ or
$ad_vw\notin\Pi$. Note that the addition to $A=ad_v$ an inner
derivation $ad_x$, $x\in\Pi$, does not change the form $\a_+$.
On the other hand the addition of $ad_w$ does change the form:
$\a_+\mapsto\a_+-\a(w)v^*$, $\a(w)\ne0$. However since $\H$ is an ideal
$dv^*=0$ and the form $\oo=d\a_+$ is preserved. As in the symplectization
we take it for the main object we have proved

 \begin{prop}\po
The set of all symplectizations is parametrized by the set
 $$
\qquad\qquad\qquad\qquad\,\quad\,\,\hfill\hfill
\{[A]\in H^1(\H,\H)\,\vert\,Aw\notin\Pi\}.
\qquad\qquad\qquad\qquad\,\quad\,\,\hfill\hfill\Box
 $$
 \end{prop}

Note that the symplectization sometimes cannot take place. Actually let
$H^1(\H,\H)=0$ (as it is in the semisimple case). Then every
derivation is inner and for $[A]=0$ we may take representative $A=0$
for which $Aw=0\in\Pi$.
\newline
%%%%%%%%%%%%%%%%%%%%%%%%%%%%%%%%%%%%%%%
\newline
{\large\bf 1.3.2.\  \
Ideals of codimension 1 in symplectic algebras}
\bigskip\smallskip

Let us consider another way to obtain a symplectic structure on a Lie
algebra $\G^{2n}=\H^{2n-1}\oplus\R^1$, $\H$ being an ideal. The Lie algebra
structure is given by the class $[A]\in H^1(\H,\H)$ of the derivation
$A=ad_v$, $v\in\R^1$. Let's assume $\H$ is equipped with a closed 2-form with
$\op{rk}\oo=2n-2$. Extend it to $\G$ by the formula $i_v\oo=\a$,
$\a\in\H^*$. Since $\left.d\oo\right\vert_\H=0$ the closeness takes the form
 $$
0=i_vd\oo=L_v\oo-di_v\oo=A\oo-d\a
 $$
(both the action and the differential live on $\H$). Let
$w\in\op{Ker}(\left.\oo\right\vert_\H)\setminus\{0\}$. The nondegeneracy of
$\oo$ on $\G$ $\a\w\oo^{n-1}\ne0$ is equivalent to
$w\notin\Pi^{2n-2}=\op{Ker}\a$. Note that changing of $A$ by an inner
derivation $ad_x$, $x\in\H$, preserves the condition $A\oo=d\a$
because under this transformation $\a\mapsto\a+i_x\oo$ and $di_x\oo=L_x\oo$
on $\H$. The condition $\a(w)\ne0$ is also preserved since $\oo(w,x)=0$ for
all $x\in\H$. However we can modify $\a$ by closed forms $\a\mapsto\a+\b$,
$d\b=0$. Let's called the described construction also the symplectization.
We have proved

 \begin{prop}\po
The set of all symplectizations of an odd-dimensional Lie algebra equipped
with maximally nondegenerate closed 2-form $\oo\in\w^2\H^*$ is parametrized
by the set
 $$
\qquad\qquad\ \hfill\hfill
\{[A]\in H^1(\H,\H),\,\a\in\H^*\,\vert\,A\oo=d\a,\,\a\w\oo^{n-1}\ne0\}.
\qquad\qquad\ \hfill\hfill\Box
 $$
 \end{prop}

Note that if $\oo$ is a differential of a contact form $\oo=d\g$,
$\g\w(d\g)^{n-1}\ne0$, then the construction just described coincides with
the symplectization from 1.3.1.

 \begin{cor}
In symplectic Lie algebra there does not exist a semisimple ideal of
codimension~1.
 \end{cor}

 \begin{proof}
In this case $[A]\in H^1(\H,\H)=0$, one can take $A=0$. Thus $d\a=0$ which
according to $H^1(\H)=0$ implies $\a=0$ and $\a\w\oo^{n-1}=0$. \qed
 \end{proof}

%%%%%%%%%%%%%%%%%%%%%%%%%%%%%%%%%%%%%%%%%%%%%%%%%%%%%%%%%%%%%%%%%%%%%%%%%%
% 1.4 %
\section{Classification in dimensions less than 5}

\hspace{13.5pt}
The first nontrivial dimension is 3 and in it almost all algebras are
contact. In what follows let $\G'$ denote the commutator of $\G$.
In the next theorem we use definition $1'$.

 \begin{th}\po
Every three-dimensional Lie algebra is contact save for the commutative
algebra and the Lie algebra $\langle x,y,z\,\vert\, [x,z]=x, [y,z]=y,
[x,y]=0 \rangle$. The contact structure is unique with the exception of
Lie algebra $sl(2)$ where we have two contact structures (up to an
isomorphism).
 \end{th}

 \begin{proof}
According to definition $1'$ a contact structure is a 2-dimensional subspace
which is not a subalgebra.
For $\op{dim}\G'=3$ we have two simple Lie algebras $so(3)$ and $sl(2)$.
For $so(3)$ the Lie multiplication is given by the vector product and all
2-planes are equivalent. For $sl(2)$ if $\Pi^2$ is not a subalgebra let
$\langle z\rangle=(\Pi^2)'$.
Note that the choice of $\Pi^2$ gives us the orientation of the space
$\G^3$. Actually for any bivector $x\w y\in\w^2\Pi^2$ it is given by the
3-vector $x\w y\w[x,y]$. Now $ad_z$ is an automorphism of $\Pi^2$. As
it is a derivation it is divergence-free, $\op{Tr}(ad_z)=0$, and we
may assume $\op{det}(ad_z)=1$. Thus $\op{Sp}(ad_z)=\{\pm i\}$ and the
claim follows. If $\op{dim}\G'=2$, $\G'$ is commutative and hence
$\op{dim}(L)=1$ for $L=\Pi^2\cap \G'$. If for a transversal to $\G'$ vector
$z$ $\left.ad_z\right\vert_{\Pi}=\op{const}\cdot\1$ the Lie algebra $\G$
is not contact. Otherwise $L^1$ may be any line in $\G'$ different from
eigenspaces of $ad_z$ and the classification follows. In the last case
$\op{dim}\G'=1$ the plane $\Pi^2$ is arbitrary transversal to $\G'$ and to
the center so that $\Pi^2$ is unique up to isomorphism. \qed
 \end{proof}

 \begin{rk}
We just showed that Lie algebra $sl(2)$ possesses two contact structures.
The Killing form $k(X,Y)$ is a nondegenerate quadratic form of the signature
$(2,1)$. So there exists a conus of isotropic vectors in $sl(2)$ and we
distinguish the contact structures subject to possibilities of $k$ being
nondegenerate on $\Pi^2=\op{Ker}\a$ or having (two) isotropic directions.
The case of one isotropic direction corresponds to $\Pi^2$ being a
subalgebra. Any 2-subspace is equally well characterized by the orthogonal
1-subspace $(\Pi^2)^\perp$ with respect to the Killing form. We may assume
that the orientation on $sl(2)$ is given by the 3-vector $X_0\w X_1\w X_2$
with
 $$
X_0=\left(
\begin{array}{cc}
1 & 0 \\ 0 & -1
\end{array}
\right),\
X_1=\left(
\begin{array}{cc}
0 & 1 \\ -1 & 0
\end{array}
\right),\
X_2=\left(
\begin{array}{cc}
0 & 1 \\ 1 & 0
\end{array}
\right).
 $$
Then we have two possibilities: the contact structure $\a$ is {\it
positive\/} $\a\w d\a>0$ which corresponds to $k(X,X)>0$ for every
$X\in(\Pi^2)^\perp\setminus\{0\}$ and the contact structure $\a$ is {\it
negative\/} $\a\w d\a<0$ which corresponds to $k(X,X)<0$.

The spherization $ST^*M^2$ of the cotangent bundle of every surface of
genus $g>0$ can be obtained as quotient of the group $\op{Sl}_2(\R)$ by
a discrete subgroup. Actually this follows from isomorphism
$\op{Sl}_2(\R)/{\{\pm1\}}\simeq ST^*L^2$,
$L^2$ being the Lobachevskii plane (see~\cite{GGPS} for details). This
isomorphism agrees with the orientation fixed above if we define the
orientation on $T^*M^2$ by the canonical symplectic form: $dv=\dfrac12
\oo^2=dp_1\w dq_1\w dp_2\w dq_2$.
Now the quotient procedure above gives us two
contact structures on $ST^*M^2$: positive $\a_+$ which is the standard
$pdq$ and negative $\a_-$ which is the connection form associated with a
metric of constant negative curvature.
 \end{rk}

Next we consider symplectic algebras of dimension~4.

 \begin{lem}\po
There does not exist a four-dimensional Lie algebra $\G=\G'$.
 \end{lem}

 \begin{proof}
Let $\G=\H\oplus R$ be a Levi decomposition, $\H$ being a semisimple
subalgebra and $R$ being a radical. We have $\op{dim}R=1$. There is an
action of $\H$ on $R$. Since $\H$ is simple the action is trivial
(otherwise the kernel is an ideal). Thus we have direct Lie summation
$\G=\H\oplus R$ and $\G'=\H'=\H$. $\Box$
 \end{proof}

The classification of 4-dimensional structures may be obtained from the
Bianchi's classification of 3-dimensional algebras and the suspension
method from 1.3.2. The proof of the following theorem uses a more direct
method. It seems rather technical but this is due to the fact that it
almost coincides with the classification of four-dimensional Lie algebras
(which may be extracted from the proof).

 \begin{th}\po
A four-dimensional Lie algebra $\G$ is symplectic iff it has the form (we
list them according to the dimension of $\G'$):

 \begin{enumerate}

\item
$\op{dim}\G'=3$. Here an algebra is symplectic iff $H^4(\G)=0$.
More precisely, the structure Maurer-Cartan equations have the form
$de_4^*=0$, $de_3^*=e_1^*\w e_2^* + e_3^*\w e_4^*$ and
  \begin{enumerate}
 \item
$de_1^*=\ll e_1^*\w e_4^*$,
$de_2^*=(1-\ll)e_2^*\w e_4^*$, $\ll\ne0,1$.
 \item
$de_1^*=\dfrac12e_1^*\w e_4^* + \ll e_2^*\w e_4^*$,
$de_2^*=-\ll e_1^*\w e_4^* +\dfrac12e_2^*\w e_4^*$.
 \item
$de_1^*=\dfrac12e_1^*\w e_4^*+\k e_2^*\w e_4^*$,
$de_2^*=\dfrac12e_2^*\w e_4^*$.
 \item
$de_1^*=-\ll e_2^*\w e_3^*$,
$de_2^*=\ll e_1^*\w e_3^*+ e_2^*\w e_4^*$, $\ll=\pm1$.
  \end{enumerate}

\item
$\op{dim}\G'=2$. In this case the structure equations are
$de_3^*=de_4^*=0$ and
  \begin{enumerate}
 \item
$de_1^*=e_1^*\we e_3^* + e_2^*\we e_4^*$,
$de_2^*= e_2^*\we e_3^* + \nu_1 e_1^*\we e_4^* + \nu_2 e_2^*\we e_4^*$.
 \item
$de_1^*=e_1^*\we e_3^*$, $de_2^*= e_3^*\we e_4^*$.
 \item
$de_1^*=e_2^*\we e_3^*$, $de_2^*= e_3^*\we e_4^*$.
 \item
$de_1^*=e_1^*\we e_3^*$, $de_2^*= -e_2^*\we e_3^*$.
 \item
$de_1^*=e_2^*\we e_3^*$, $de_2^*= -e_1^*\we e_3^*$.
  \end{enumerate}

\item $\op{dim}\G'=1$. Here $de_2^*=de_3^*=de_4^*=0$ and
  \begin{enumerate}
 \item
$de_1^*= e_1^*\we e_2^*$.
 \item
$de_1^*= e_2^*\we e_3^*$.
  \end{enumerate}

\item
Lie algebra $\G$ is commutative, $\G'=0$.

 \end{enumerate}

In addition in cases 1 and 2$\text{(i)}$ the symplectic structure can be
chosen exact and the set of such structures is nonempty Zarissky-open
in $\op{Im}\left.d\right|_{C^1}\subset C^2$, in particular everywhere dense.
In other cases the symplectic form is unique (in 2$\text{(ii)}$ up to the
sign) and is not exact. In cases 2$\text{(ii)}$-$\text{(iii)}$ the structure
has the form $\oo=e_1^*\we e_3^* + e_2^*\we e_4^*$, and in cases
2$\text{(iv)}$-$\text{(v)}$, 3 and 4 it is
$\oo=e_1^*\we e_2^* + e_3^*\we e_4^*$.
 \end{th}

 \begin{proof}
${\bf 1}^\circ$.
\underline{$\op{dim}\G'=3$}.
Thanks to lemma~1 there exists a covector $e_4^*$ such that $de_4^*=0$.
First consider the case all exact 2-forms $df^*$ are degenerate, i.e.
decomposable. Then there exists a 3-dimensional subspace
$\Pi^3\subset\G^*$ such that $dC^1\subset\we^2\Pi^3\subset C^2$.

Consider the case $e_4^*\in\Pi^3$. One may complete $e_4^*$ to a basis
$e_i^*$ such that $\Pi^3=\la e_2^*, e_3^*, e_4^*\ra$ and $de_2^*=f^*\we
e_4^*$, with $f^*\in\la e_2^*, e_3^*\ra$. Up to a change of basis
there are two possibilities: $f^*=e_2^*$, $de_3^*\notin\la e_2^*,
e_3^*\ra\we e_4^*$ or $de_3^*=h^*\we e_4^*$, $h^*\in\la e_2^*, e_3^*\ra$.
In the first case one chooses a covector $e_1^*$ satisfying
$de_1^*=e_3^*\we e_4^*$. Then $de_3^*=a e_2^*\we e_3^* + v^*\we
e_4^*$, $a\ne0$, and the equation $d^2e_1^*=0$ does not hold.
In the second case the covector $e_1^*$ can be chosen to satisfy the
equality $de_1^*=e_2^*\we e_3^*$. Let us consider a linear automorphism
of the plane $\la e_2^*, e_3^*\ra$, given by the formulas: $e_2^*\mapsto f^*$,
$e_3^*\mapsto h^*$. The equation $d^2e_1^*=0$ implies its trace vanishes.
Moreover eigenvalues are nonzero and are complex conjugated if not real.
Thus after a transformation $de_2^*=e_2^*\we e_4^*$ and
$de_3^*=-e_3^*\we e_4^*$, or $de_2^*= e_3^*\we e_4^*$ and
$de_3^*= -e_2^*\we e_4^*$. It's easy to see that in every case the Lie
algebra is not symplectic and moreover $H^4(\G)\ne0$.

Let now $e_4^*\notin\Pi^3$. We may assume that $\Pi^3$ is generated by
covectors $e_1^*$, $e_2^*$ and $e_3^*$, i.e. there is a direct Lie
decomposition $\G=\H\oplus\R^1$, $\H$ being a 3-dimensional Lie algebra.
We may choose the basis in $\Pi^3$ in such a manner that the symplectic
form $\oo=e_1^*\we e_2^*+ e_3^*\we e_4^*$. Then $d\oo=0$ implies
$de_3^*=0$. Thus $\op{dim}\G'<3$. According to Bianchi's classification
there are exactly two 3-dimensional algebras with the commutator
$\H=\H'=\G'$: $\H=sl(2)$ and $\H=so(3)$ and in every case
$H^4(\G)=H^3(\H)=\R$.

Let us complete the classification of symplectic algebras $\G$ of the type
$\H\oplus\R$. We have:
$de_1^*=A_1e_1^*\we e_2^*+ B_1e_1^*\we e_3^*+ C_1e_2^*\we e_3^*$,
$de_2^*=B_2e_1^*\we e_3^*+ C_2e_2^*\we e_3^*$, $de_3^*=de_4^*=0$. The
condition $d\oo=0$ gives $B_1+C_2=0$ and the conditions $d^2e_i^*=0$
are equivalent to the equalities $A_1B_1=A_1B_2=0$. Assume first
$A_1=0$. Considering the canonical forms of linear automorphism of the
plane $\la e_1^*, e_2^*\ra$, given by the formula $e_i^*\mapsto
de_i^*(\cdot,e_3^*)$, we get the Lie algebras 2(iv), 2(v) and 3(ii).
The other case $A_1\ne0$ leads to the Lie algebra 3(i).
Let us call further the considered case 3-dimensional.

Now let us consider the possibility of nondegenerate form $df^*$,
$f^*\in\G^*$. It is clear that if Lie algebra $\G$ possesses an exact
symplectic form $df^*$, then for a Zarissky-open set of covectors
$f^*\in\G^*$ the 2-form $df^*$ is symplectic.

We may assume that there exists a basis $e_i^*$ such that
$de_3^*=e_1^*\w e_2^*+e_3^*\w e_4^*$. Actually this is equivalent to the
existence of a covector $X\notin\la e_4^*\ra$ and a number $\ll\ne0$ such
that $\sigma=dX-\ll X\w e_4^*$ is a decomposable 2-form, i.e.
$0=\sigma^2=(dX)^2-\ll dX\w X\w e_4^*$. For almost every $X$ the form
$(dX)^2$ is nondegenerate. So the pair $(X,\ll)$ does exist iff there exists
a covector $X$ such that $dX\w X\w e_4^*\ne0$. Let us suppose this 4-form is
zero for every $X$. Let us consider some generic covector $X=e_1^*\notin\la
e_4^*\ra$. Then by our arguments there are covectors $e_2^*$, $e_3^*$ with
the property $de_1^*=e_1^*\w e_2^*+ e_3^*\w e_4^*$. The differentials of
the other basis covectors must have the form $de_2^*=e_2^*\w X_1+e_4^*\w
X_2$, $de_3^*=e_3^*\w X_3+ e_4^*\w X_4$, $X_i\in\la e_1^*,e_2^*,e_3^*\ra$.
From the equation $d^2e_1^*=0$ we have: $X_1=\mu e_1^*$, $X_2=\nu e_1^*+\k
e_3^*$, $X_3=\k e_1^*+e_2^*$. Thus $de_2^*=-\mu e_1^*\w e_2^*-\nu e_1^*\w
e_4^*-\k e_3^*\w e_4^*$ and from $d^2e_2^*=0$ we get $\nu=0$, i.e. as
$d(e_1^*\w e_2^*)\ne0$ $de_2^*=-\mu de_1^*$ which contradicts the
assumption $\op{dim}\G'=3$.

Now let us consider the general form for the differentials (we may assume
$A_2=0$):
 \begin{eqnarray*}
de_1^*= & \!\! A_1e_1^*\w e_2^*+ & \!\!\!\!
B_1e_1^*\w e_3^*+ C_1e_1^*\w e_4^*+ D_1e_2^*\w e_3^*+ E_1e_2^*\w e_4^*+
F_1e_3^*\w e_4^*,\\
de_2^*= & \!\!                   & \!\!\!\!
B_2e_1^*\w e_3^*+ C_2e_1^*\w e_4^*+ D_2e_2^*\w e_3^*+ E_2e_2^*\w e_4^*+
F_2e_3^*\w e_4^*.
 \end{eqnarray*}

The condition $d^2e_3^*=0$ implies $F_1=F_2=0$, $C_1+E_2=1$, $B_1+D_2=0$.
From $d^2e_2^*=0$ we have $B_2A_1=C_2A_1=0$, and
 $$
-B_2C_1-B_2+C_2B_1-D_2C_2+E_2B_2=0,\ -B_2E_1+C_2D_1-D_2=0.  \eqno (\dagger)
 $$
If $A_1\ne0$ then $B_2=C_2=D_2=0$ and $d^2e_1^*=0$ gives $E_2=0$ which
contradicts $de_2^*\ne0$. Thus $A_1=0$. Consider the automorphism of the
plane $\la e_1^*,e_2^*\ra$ given by the formula $e_i^*\mapsto
de_i^*(\cdot,e_3)$ and consider its canonical forms (according to the
equations above it is traceless).

a). $B_1=-D_2=\ll$, $D_1=B_2=0$. Then $(\dagger)$ implies $\ll=0$.

b). $B_1=B_2=D_2=0$, $D_1=\ll$. If $\ll\ne0$ then $(\dagger)$ implies
$C_2=0$. Now the equation $d^2e_1^*=0$ implies $E_2=0$ which is impossible
since $de_2^*\ne0$.

c). $B_1=D_2=0$, $B_2=-D_1=\ll$. Let us suppose $\ll\ne0$. Then $(\dagger)$
implies $E_1=-C_2$, $C_1=0$, $E_2=1$. Making the transformation
$e_3^*=e_3^*-\dfrac{E_1}\ll e_4^*$ we obtain case 1$\text{(iv)}$.

Otherwise we have $B_i=D_i=0$. Hence we can consider
the automorphism of the plane $\la e_1^*,e_2^*\ra$ given by the formula
$e_i^*\mapsto de_i^*(\cdot,e_3)$ and its canonical forms. Since the trace
of this automorphism is $C_1+E_2=1$ we obtain all cases
1$\text{(i)}$-$\text{(iii)}$.

${\bf 2}^\circ$.
\underline{$\op{dim}\G'=2$}.
For some basis we have $de_3^*=de_4^*=0$ and differentials
$de_1^*$ and $de_2^*$ are linear independent and have a general form:
 \begin{eqnarray*}
de_1^*
\!\!\!\!&=&\!\!\!\!
          A_1e_1^*\we e_2^* + B_1e_1^*\we e_3^* + C_1e_1^*\we e_4^*
        + D_1e_2^*\we e_3^* + E_1e_2^*\we e_4^* + F_1e_3^*\we e_4^*,\\
de_2^*
\!\!\!\!&=&\!\!\!\!
 \phantom{A_2e_1^*\we e_2^* + \ }
                              B_2e_1^*\we e_3^* + C_2e_1^*\we e_4^*
        + D_2e_2^*\we e_3^* + E_2e_2^*\we e_4^* + F_2e_3^*\we e_4^*.
 \end{eqnarray*}

From the conditions $d^2e_i^*=0$ we have:
 $$
A_1F_1-B_1E_1+C_1D_1-E_2D_1+E_1D_2=0,
 $$
 $$
-A_1F_2-D_1C_2+E_1B_2=0,\  -B_2C_1+B_1C_2-C_2D_2+B_2E_2=0, \eqno (\dd)
 $$
 $$
A_1B_2=A_1C_2=A_1D_2=A_1E_2=0.
 $$

If $A_1\ne0$ then $B_2=C_2=D_2=E_2=F_2=0$ and $de_2^*=0$, i.e.
$\op{dim}\G'=1$. Thus $A_1=0$.

Let us consider an automorphism of the plane
$\la e_1^*,e_2^*\ra$ given by the formula
$e_i^*\mapsto de_i^*(\cdot,e_3)=B_ie_1^*+ D_ie_2^*
(\op{mod}e_4^*)$, $i=1,2$. Further we consider all possible canonical
forms of this automorphism. The repeated phrase "we may assume" means
"there exists a coordinate transformation such that".

1. The eigenvalues $\ll_1\ne\ll_2$. Let, say, $\ll_1\ne0$. We may assume
$B_1=\ll_1$, $B_2=0$, $D_1=0$, $D_2=\ll_2$. From the conditions $(\dd)$
we have: $C_2=E_1=0$. We may assume $C_1=F_1=0$.

1.1. $E_2\ne0$. By a transformation is reduced to case
2(i), $\nu_1=0$, $\nu_2=1$.

1.2. $E_2=0$, $\ll_2\ne0$. By a transformation is reduced to
the 3-dimensional case considered above.

1.3. $E_2=\ll_2=0$. In this case $F_2\ne0$ and we get case 2(ii).

2. Jordan box with nonzero eigenvalue: $B_1=1$, $B_2=0$,
$D_1=\m$, $D_2=1$. From the conditions $(\dd)$ we have:
$C_2\m=0$, $(C_1-E_2)\m=0$.

2.1. $\m=0$. We may assume $F_1=F_2=0$. A transformation in the plane
$\la e_3^*, e_4^*\ra$ leads to $C_1=0$.

2.1.1. $E_1\ne0$. In this case we have equations 2(i).

2.1.2. $E_1=0$, $C_2\ne0$. A transformation leads to case 2(i).

2.1.3. $E_1=C_2=0$, $E_2\ne0$. By a transformations is reduced to 2(i),
$\nu_1=0$, $\nu_2=1$.

2.1.4. $C_2=E_1=E_2=0$. We get the 3-dimensional case considered above.

2.2. $\m\ne0$. Hence $C_2=0$, $C_1=E_2$. A transformation in the space
$\la e_2^*, e_3^*, e_4^*\ra$ allows to assume $de_2^*=e_2^*\we e_3^*$.
With this change the coefficients of the decomposition of
$de_1^*$ become arbitrary save for the conditions $A_1=0$ and
$E_1-B_1E_1+C_1D_1=0$.

2.2.1. $B_1=1$, $C_1D_1=0$. We may assume $F_1=0$.

2.2.1.1. $C_1=0$, $E_1\ne0$. We get a special case of 2(i).

2.2.1.2. $C_1=0$, $E_1=0$. We get the 3-dimensional case considered above.

2.2.1.3. $D_1=0$, $E_1\ne0$. A transformation leads to case 2(i).

2.2.1.4. $D_1=0$, $E_1=0$, $C_1\ne0$. A transformation leads to 2(i),
$\nu_1=0$, $\nu_2=1$.

2.2.1.5. $D_1=0$, $E_1=0$, $C_1=0$. We get the 3-dimensional case considered
above.

2.2.2. $B_1\ne1$, $B_1\ne0$, $E_1=\dfrac{C_1D_1}{B_1-1}$. We may assume
$F_1=0$.

2.2.2.1. $C_1D_1\ne0$. We may assume $C_1=B_1-1$. We get:
$de_1^*=B_1e_1^*\we e_3^* + (B_1-1)e_1^*\we e_4^* +D_1e_2^*\we e_3^*
+D_1e_2^*\we e_4^* = e_1^*\we e_3^* + D_1e_2^*\we (e_3^* + e_4^*) +
(B_1-1)e_1^*\we(e_3^*+e_4^*)$, and after a transformation:
$de_1^*= e_1^*\we e_3^* + e_2^*\we e_4^* + \nu e_1^*\we e_4^*$.
May be lead to equations 2(i).

2.2.2.2. $C_1\ne0$, $D_1=0$. A transformation leads to 2(i), $\nu_1=0$,
$\nu_2=1$.

2.2.2.3. $C_1=0$. We get the 3-dimensional case considered above.

2.2.3. $B_1=0$, $E_1=-C_1D_1$.

2.2.3.1. $C_1D_1\ne0$. We may assume $C_1=1$, $F_1=0$. We get
$de_1^*=e_1^*\we e_4^* + D_1e_2^*\we e_3^* - D_1e_2^*\we e_4^*$. Making
the change $D_1(e_3^*-e_4^*)\mapsto e_3^*$ and after
$e_3^*\leftrightarrow e_4^*$ we get equations 2(i).

2.2.3.2. $C_1\ne0$, $D_1=0$. We may assume $C_1=1$, $F_1=0$.
A transformation leads to 2(i), $\nu_1=0$, $\nu_2=1$.

2.2.3.3. $C_1=0$. A transformation leads to 2(ii).

3. Jordan box with zero eigenvalue: $B_1=0$, $B_2=0$,
$D_1=\m$, $D_2=0$. From the conditions $(\dd)$ we have:
$C_2\m=0$, $(C_1-E_2)\m=0$.

3.1. $\m=0$. We have: $de_i^*=v_i^*\we e_4^*$, $v_i^*\in \la e_1^*, e_2^*,
e_3^*\ra$, $i=1,2$.

3.1.1. $e_3^*\in \la v_1^*, v_2^*\ra$.
In this case there exists a basis such that $de_1^*=v^*\we e_4^*$,
$de_2^*=e_3^*\we e_4^*$, $v^*\in\la e_1^*, e_2^*\ra$. Up to changes there
are two possibilities: $v^*=e_1^*$ and $v^*=e_2^*$. They lead to
equations 2(ii) and 2(iii) correspondingly if use the transformation
$e_4^*\mapsto e_3^*$, $e_3^*\mapsto -e_4^*$.

3.1.2. $v_i^*\in\la e_1^*, e_2^*\ra$. We get the 3-dimensional case
considered above.

3.2. $\m\ne0$. Then $C_1=E_2$, $C_2=0$.

3.2.1. $C_1\ne0$. We may assume $F_2=0$. After the transformation
$e_3^*\leftrightarrow e_4^*$ we get the case considered in 2.2 above.

3.2.2. $C_1=0$. Then $F_2\ne0$. A transformation leads to $F_1=0$, and
another one leads to $E_1=0$. We get equations 2(iii).

4. Pure imaginary conjugated roots: $B_1=0$, $B_2=-1$, $D_1=1$, $D_2=0$.
From the equations $(\dd)$ we get: $C_1=E_2$, $C_2=-E_1$.

4.1. $C_1\ne0$. We may assume $C_1=1$. A transformation leads to
$F_1=0$, and another leads to $C_2=0$. We get the equations:
$de_1^*=e_1^*\we e_4^* + e_2^*\we e_3^*$,
$de_2^*= -e_1^*\we e_3^* + e_2^*\we e_4^* + F_2 e_3^*\we e_4^*$. We may
assume $F_2=0$. The transformation $e_3^*\leftrightarrow e_4^*$ leads
to a special case of 2(i).

4.2. $C_1=0$. We may assume $F_1=0$, $C_2=0$. We get case 2(v).

5. Conjugated roots not belonging to coordinate axis:
$B_1=1$, $B_2=-a$, $D_1=a$, $D_2=1$, $a\ne0$. From the conditions $(\dd)$ we
have: $C_1=E_2$, $C_2=-E_1$. We may assume $F_1=0$. Since $a\ne0$, we may
assume $C_2=0$.

5.1. $C_1\ne0$. We may assume $C_1=1$. After the transformation
$e_3^*+e_4^*\mapsto e_4^*$ we get the equations:
$de_1^*=e_1^*\we e_4^* + ae_2^*\we e_3^*$,
$de_2^*= -ae_1^*\we e_3^* + e_2^*\we e_4^* + F_2 e_3^*\we e_4^*$. We may
assume $a=1$, $F_2=0$. After the transformation $e_3^*\leftrightarrow
e_4^*$ we get a special case of equations 3(i).

5.2. $C_1=0$. The equations take the form: $de_1^*=e_1^*\we e_3^* +
ae_2^*\we e_3^*$, $de_2^*= -ae_1^*\we e_3^* + e_2^*\we e_3^* +
F_2 e_3^*\we e_4^*$. One easily checks that in this case
the Lie algebra is not symplectic.

${\bf 3}^\circ$.
\underline{$\op{dim}\G'=1$}.
In this case the only nonzero differential is degenerate.
Actually if not then there exists a basis such that
$de_1^*=e_1^*\we e_2^* + e_3^*\we e_4^*$, $de_2^*=de_3^*=de_4^*=0$,
from where $d^2e_1^*\ne0$. Thus,
the differential $de_1^*$ is decomposable and we get
the 3-dimensional case considered above.

${\bf 4}^\circ$.
\underline{$\op{dim}\G'=0$}. This is clearly the commutative case.

The statement on the uniqueness of the structures at the end of the theorem
could be now easily checked. \qed
 \end{proof}

Now when we have an exact symplectic algebra there may exist nonequivalent
symplectic structures. Modulo closed 1-forms the set of exact structures
is parametrized by covectors whose differentials are nondegenerate. These
covectors among all covectors are organized in orbits of the coadjoint
action.
\newline
\newline

 \begin{lem}\po
If a covector $f^*\in\G^*$ is such that $df^*$ is nondegenerate then the
orbit of coadjoint action through it is open. In other words for any
sufficiently close covector ${\hat f}^*$ there exists an automorphism of the
Lie algebra such that it sends $f^*$ to ${\hat f}^*$.
 \end{lem}

 \begin{proof}
The tangent space to the orbit through $f^*$ is $T_{f^*}{\cal O}=\la
ad_X^*(f^*)\,\vert\,X\in\G \ra$. Its annulator is the linear space of
$Y\in\G$ such that $\la ad_X^*(f^*),Y\ra=0$ for every $X\in\G$. This means
that $\la f^*, [X,Y]\ra=0$ or $df^*(X,Y)=0$ for all $X$ which yields $Y=0$
and $T_{f^*}{\cal O}=\G^*$. \qed
 \end{proof}

%%%%%%%%%%%%%%%%%%%%%%%%%%%%%%%%%%%%%%%%%%%%%%%%%%%%%%%%%%%%%%%%%%%%%%%%%%
% 2 %
 \chapter{Elliptic Monge-Amp\'ere equations}

\hspace{13.5pt}
According to the paper~\cite{L} Monge-Amp\'ere equations may be considered
as half-dimensional effective forms on the elements of some contact
distribution (it's better to call such an object the generalized
Monge-Amp\'ere equation). The pointwise classification
(i.e. with respect to the linear
symplectic group at a point) of such equations was given in~\cite{LRC}.
We consider only the case of Monge-Amp\'ere with two independent variables.
This means that we consider a four-dimensional symplectic manifold
$(M^4,\oo)$ and a 2-form $\te\in\O^2(M)$, $\te\w\oo=0$.

%%%%%%%%%%%%%%%%%%%%%%%%%%%%%%%%%%%%%%%%%%%%%%%%%%%%%%%%%%%%%%%%%%%%%%%%%%
% 2.1 %
\section{Classification of elliptic equations depending on two variables}

\hspace{13.5pt}
For the elliptic Monge-Amp\'ere equations with two variables there exists
an alternative description. Let normalize $\te$ by the condition
$\op{Pf}(\te)=1$ and define an automorphism $j$ by the formula $i_j\oo=\te$
where $i_\bullet$ is the substitution. Then $j$ is an almost complex
structure. The Monge-Amp\'ere equation $(\oo,j)$ is called nondegenerate if
$j$ is of general position (at least it is nonintegrable and the equation
is not of divergence type). For these equations the paper~\cite{K} contains
the following classificational result:

 \begin{th}\po
A nondegenerate elliptic Monge-Amp\`ere equation $(\oo,j)$ canonically
determines an $\{e\}$-structure, i.e. the field of basis frames
$(P_1,P_2,Q_1,Q_2)$. This structure is a complete invariant, i.e. two
nondegenerate elliptic Monge-Amp\`ere equations are isomorphic if and only
if the corresponding $\{e\}$-structu\-res are. The classifying
$\{e\}$-structure satisfies the following relations which completely
determine it (we use the dual basis in the cotangent space):
 $$
\oo=P_1^*\w Q_1^*+P_2^*\w Q_2^*, \quad \te=P_1^*\w Q_2^*-P_2^*\w Q_1^*,
 $$
 $$
j=P_1^*\ot P_2-P_2^*\ot P_1- Q_1^*\ot Q_2+ Q_2^*\ot Q_1,
 $$
 $$
N_j=-P_1^*\w Q_1^*\ot P_2+ P_1^*\w Q_2^*\ot P_1- P_2^*\w Q_1^*\ot P_1-
P_2^*\w Q_2^*\ot P_2.
 $$
 \end{th}

Here $N_j=[j,j]$ is the Nijenhuis self-bracket of the almost complex
structure $j$. Note that since for $\{e\}$-structures the equivalence
problem is solved (see~\cite{S}) this theorem serves as equivalence
criterion for two-dimensional elliptic nondegenerate Monge-Amp\'ere
equations.

 \begin{cor}
Every symmetry of a Monge-Amp\'ere equation is a symmetry for its
$\{e\}$-structure invariant and vice versa. \qed
 \end{cor}

Thus we are given a tool for constructing symmetries for the equations of
the described type. Moreover this is the key idea for considering the
structures of the next section.

%%%%%%%%%%%%%%%%%%%%%%%%%%%%%%%%%%%%%%%%%%%%%%%%%%%%%%%%%%%%%%%%%%%%%%%%%%
% 2.2 %
\section{Lie group action and invariant equations}

\hspace{13.5pt}
It is natural to consider equivariant equations, i.e. the equations with
transitive action of some Lie group. We assume that the elements of this
group depend on the first derivatives of the solutions.
Let us call equations of these type {\it invariant equations\/}.

In our model we must permit a Lie group action on the symplectic manifold
$(M^4,\oo)$ which preserves the structures $\oo$ and $j$. Since the action
is transitive our manifold is homogeneous. Let us note that according to
the canonicity of formulas from theorem~4 this is equivalent to the
invariance of the classifying $\{e\}$-structure. Thus the structural
functions $c_{ij}^k$ for $\{e\}$-structure $e_i$, $[e_i,e_j]=c_{ij}^ke_k$,
are constant and our manifold becomes a Lie group.
As in chapter~1 we may assume everything to live on the
corresponding Lie algebra. For example the Nijenhuis tensor
 $$
N_j(\x,\e)=[j\x,j\e]-j[j\x,\e]-j[\x,j\e]-[\x,\e]
 $$
can be calculated by means of the Lie algebra commutators.

Theorem~4 now also takes place and the construction of the
$\{e\}$-structure from~\cite{K} may be re-written for the invariant
situation in the interior Lie algebra terms.

Note that the $\{e\}$-structure appeared is not arbitrary. There are two
differential conditions on it. First the form $\oo$ is closed.
Setting $e_1=P_1$, $e_2=P_2$, $e_3=Q_1$, $e_4=Q_2$ from the structural
equations above wee obtain:
 $$
c^1_{12}-c^3_{23}+c^4_{13}=c^2_{12}-c^3_{24}+c^4_{14}=
c^1_{14}-c^2_{13}+c^3_{34}=c^1_{24}-c^2_{23}+c^4_{34}=0.
 $$
The second condition is connected with the almost complex structure. We may
define this structure $j$ by means of formulas of theorem~4 and then we
compute the Nijenhuis tensor. The condition is that it coincides with the
tensor $N_j$ given in theorem~4. If these two conditions hold true we may
recover the Monge-Amp\'ere equation.

%%%%%%%%%%%%%%%%%%%%%%%%%%%%%%%%%%%%%%%%%%%%%%%%%%%%%%%%%%%%%%%%%%%%%%%%%%
% 2.3 %
\section{Example of the recovering a Monge-\protect\\
Amp\'ere equation by its $\{e\}$-structure}

\hspace{13.5pt}
Let two necessary conditions discussed in 2.2 be satisfied.
Define coordinates in a neighborhood of the unity basing on the exponential
mapping $\op{exp}: {\cal G}\to G$. Campbell-Hausdorff formula shows how
to write left-invariant vector fields in these coordinates (\cite{J},
\cite{SL}). This allows us to write down the generalized Monge-Amp\'ere
equation. To obtain the ordinary one (\cite{L}) we need to fix a Lagrange
submanifold $L^2\subset M^4$ (or a local cotangent bundle with $\te$-trivial
Lagrange fibers $L^2$ in order to obtain a quazilinear equation), identify
its neighborhood with cotangent bundle $UL\simeq T^*L$ and substitute the
expression $p=\p u/ \p q$ into the equation $\te(p,q)=0$ in canonical
coordinates $(p,q)$. The last operation may result in different
Monge-Amp\'ere for different choices of $L$ (however for nonequivalent
admissible $\{e\}$-structures or for nonequivalent generalized
Monge-Amp\'ere equations all possible as representative ordinary
Monge-Amp\'eres are different).

Let's demonstrate the scheme for a nilpotent Lie algebra (we may call the
obtained equation nilpotent). A nilpotent Lie algebra on $\R^4$ is
isomorphic to one of the cases: commutative, 3-dimensional or 4-dimensional.
Let us consider the last case. It is determined by the relations:
$[e_2,e_3]=e_1$, $[e_3,e_4]=e_2$ and $[e_i,e_j]=0$ for all the others
$i<j$. Let us consider the left-invariant basis of the Lie algebra of the
form: $P_1=e_1$, $P_2=e_2$, $Q_1=e_3$, $Q_2=e_4$. This basis satisfies two
necessary conditions from the end of 2.2. Write it in exponential
coordinates. Note that for 4-dimensional nilpotent case the series in the
Campbell-Hausdorff formula terminates on the third term. Thus the
left-invariant vector field through a point $x\in{\cal G}$
taking the value $\e\in{\cal G}$ at zero has the form
($L_x$ denotes the left shift):
 $$
(L_x)_*\e=\e+\frac12[x,\e]+\frac1{12}[x,[x,\e]].
 $$
Extending left-invariantly the basis $e_i=\p_i=\dfrac\p{\p x_i}$ we get
the expression for the basis in coordinates $x$:
 \begin{eqnarray*}
&
P_1=\p_1,\qquad\quad\ \,
&
Q_1=\p_3+(\dfrac12x_2+\dfrac1{12}x_3x_4)\p_1-\dfrac12x_4\p_2,\\
&
P_2=\p_2-\dfrac12x_3\p_1,
&
Q_2=\p_4-\dfrac1{12}(x_3)^2\p_1+\dfrac12x_3\p_2.
 \end{eqnarray*}

For the dual basis we have the expression:
 \begin{eqnarray*}
&
P_1^*=dx_1+\dfrac12x_3dx_2-(\dfrac12x_2-\dfrac16x_3x_4)dx_3-
\dfrac16(x_3)^2dx_4,\
&
Q_1^*=dx_3,\\
&
P_2^*=dx_2+\dfrac12x_4dx_3-\dfrac12x_3dx_4,
\qquad\qquad\qquad\qquad\qquad\ \,
&
Q_2^*=dx_4.
 \end{eqnarray*}

Note that the last expressions may be obtained from the Maurer-Cartan
equations $de^*_1=e_3^*\w e_2^*$, $de_2^*= e_4^*\w e_3^*$,
$de_3^*=de_4^*=0$. Further:
 \begin{eqnarray*}
\oo &=& P_1^*\wedge Q_1^* + P_2^*\wedge Q_2^* \\
    &=& dx_1\we dx_3+dx_2\we dx_4+\frac12x_3dx_2\we dx_3\\
&& +(\dfrac12x_4+\dfrac16(x_3)^2)dx_3\we dx_4.
 \end{eqnarray*}
 \begin{eqnarray*}
\te &=& P_1^*\we Q_2^* - P_2^*\we Q_1^*  \\
 &=& dx_1\we dx_4 - dx_2\we dx_3+\frac12x_3dx_2\we dx_4\\
&& +(-\frac12x_2-\frac12x_3+\frac16x_3x_4)dx_3\we dx_4.
 \end{eqnarray*}

Let change coordinates:

 \begin{eqnarray*}
&
p_1=x_1+\dfrac12x_4+\dfrac14x_3^2+\dfrac12x_2x_3-\dfrac12x_4^2-
\dfrac16x_3^2x_4, \
&
q_1=x_3, \\
& p_2=x_2+\dfrac12x_3-\dfrac12x_3x_4,
\qquad\qquad\qquad\qquad\qquad
&
q_2=x_4.
 \end{eqnarray*}

This choice is equivalent to the choice of the general solution
$L^2=\{q_1=c_1; q_2=c_2\}$, i.e. such a Lagrangian manifold
that $\left.\te\right|_L\equiv0$. In new coordinates we have:
 $$
\oo=dp_1\we dq_1+ dp_2\we dq_2.
 $$

 $$
\te=dp_1\we dq_2 -dp_2\we dq_1-p_2dq_1\we dq_2.
 $$

Substitute the expressions $p_1=\p u/\p q_1$, $p_2=\p u/\p q_2$ to the
equation $\te=0$ and we get the Monge-Amp\'ere equation which we looked for:

 $$
\triangle u=\dfrac{\p u}{\p q_2}.
 $$

%%%%%%%%%%%%%%%%%%%%%%%%%%%%%%%%%%%%%%%%%%%%%%%%%%%%%%%%%%%%%%%%%%%%%%%%%%
% 2.4 %
\section{Scheme of normal forms for Monge-\protect\\
Amp\'eres on two-dimensional surfaces}

\hspace{13.5pt}
Now since theorem~4 gives an equivalence criterion for nondegenerate
elliptic Monge-Amp\'ere equations with two variables it's possible to
determine normal forms of these equations. We suppose that the situation is
equivariant and so we have a Lie group with invariant symplectic and almost
complex structures. As usual we talk instead of Lie algebras.

We may use normal forms of symplectic structures on Lie algebras as in 1.4.
Then we add an almost complex structure $j$ to the pair $(\G,\oo)$ with $\G$
a Lie algebra and $\oo$ a symplectic form. Now this form is not arbitrary.
It must satisfy the condition $\oo(X,jX)=0$. Denote by $\J$ the set of all
complex structures on a linear space $V$ and by $\J^\pm(V,\oo)$ the set
$\{j\in\J(V)\,\vert\,\oo(jX,jY)=\pm\oo(X,Y)\}$.
The set $\J^+(V,\oo)$ is not even connected (contrary to
the set $\J^+(V,\oo)\cap\{j\in\J\,\vert\,\oo(X,jX)>0,\,X\ne0\}$
which as well-known in symplectic geometry is contractible).
We are interested in
the set $\J^-(V^4,\oo)$. This manifold is diffeomorphic to $S^2\t \R^2$.
Actually instead of $j$ we may consider 2-forms $\te(X,Y)=\oo(jX,Y)$ which
satisfy the conditions $\te\w\oo=0$, $\la\te,\te\ra=1$. Here
$\la\te_1,\te_2\ra=\dfrac{\te_1\w\te_2}{\oo\w\oo}$ is the metric associated
with the Pfaffian. This metric has type $(3;3)$. The orthogonal complement to
$\oo$ has type $(2;3)$ and now the equality $\J^-(\G^4,\oo)\simeq S^2\t\R^2$
is clear. Note that actually this gives an obstruction to existence of a
nonvanishing generalized elliptic Monge-Amp\'ere equation on a symplectic
manifold $(M^4,\oo)$, this is the first obstruction for constructing a
section to the bundle $\J^-(T_xM^4;\oo_x)\mapsto x$ and it lies in
$H^3(M^4;\Z)$. However in the case of Lie groups this obstruction is
always zero. Actually there always exists a left-invariant section.

Thus to the normal symplectic forms from 1.4 we can add the structures from
$\J^-(\G^4;\oo)$ and this would gives us through the recovering procedure
of 2.3 the normal forms of elliptic Monge-Amp\'ere equations with two
variables.

%\newpage
%
%
%
%%%%%%%%%%%%%%%%%%%%%%%%%%%%%%%%%%%%%%%%%%%%%%%%%%%%%%%%%%%%%%%%%%%%%%%%%%
%
%%%%%%%%%%%%%%%%%%%%%%%%%%%%%%%%%%%%%%%%%%%%%%%%%%%%%%%%%%%%%%%%%%%%%%%%%%
%

%

\bigskip
\bigskip
\bigskip

\ {\hbox to 12.5cm{ \hrulefill }}

\bigskip
\bigskip
\bigskip

{\it \hspace{-19pt} Address:}
{\footnotesize
 \begin{itemize}
  \item
P. Box 546, 119618, Moscow, Russia
  \item
Chair of System Analysis, Moscow State Technological University
\linebreak
{\rm n. a.} Baumann, Moscow, Russia
 \end{itemize}
}

{\it \hspace{-19pt} E-mail:} \quad
{\footnotesize
lychagin\verb"@"glas.apc.org or borkru\verb"@"difgeo.math.msu.su
}

%%%%%%%%%%%%%%%%%%%%%%%%%%%%%%%%%%%%%%%%%%%%%%%%%%%%%%%%%%%%%%%%%%%%%%%%%%
\end{document}